\documentclass[a4paper,12pt]{article}
\usepackage{latexsym}

\usepackage{amssymb} % \mathbb
\usepackage{amsmath} % \mathbb
\usepackage{theorem}

\def\P{{\mathbf{P}}}% \P == \mathbb{P}
\def\Z{{\mathbb{Z}}}% \Z == \mathbb{Z}
\DeclareMathOperator{\rank}{rank}

\DeclareMathOperator{\Sing}{Sing}
\DeclareMathOperator{\depth}{depth}

\DeclareMathOperator{\Supp}{Supp}

\DeclareMathOperator{\codim}{codim}

\DeclareMathOperator{\Hom}{Hom}
\DeclareMathOperator{\Ext}{Ext}
\DeclareMathOperator{\Proj}{Proj}

\DeclareMathOperator{\pd}{pd}

\DeclareMathOperator{\Der}{Der}
\def\K{{\mathbb{K}}}% \P == \mathbb{P}

\newcommand{\owari}{\hfill$\square$}

\theoremstyle{break}
\newtheorem{theorem}{Theorem}[section]
\newtheorem{prop}[theorem]{Proposition}
\newtheorem{cor}[theorem]{Corollary}
\newtheorem{lemma}[theorem]{Lemma}
\newtheorem{define}{Definition}[section]
\newtheorem{rem}{Remark}[section]

 \title{Splitting criterion for reflexive sheaves}
 \author{TAKURO ABE \and MASAHIKO YOSHINAGA}
 \date{\today}
\begin{document}
 \maketitle

 \begin{abstract}
The purpose of this paper is to study the structure of  
reflexive sheaves over projective spaces through 
hyperplane sections. 
We give a criterion 
for a reflexive sheaf to 
split into a direct sum of line bundles. 
An application to the theory of free hyperplane arrangements 
is also given.
 \end{abstract}

 \setcounter{section}{-1}
 \section{Main Theorem}
 \label{intro}

Vector bundles over the projective space 
$\P^n_\K$ are one of the main subjects in both (algebraic) geometry
and commutative algebra. The most fundamental result in this area is
the theorem due to Grothendieck which asserts that any
holomorphic vector bundle over $\P^1_\K$ splits into a direct
sum of line bundles.
When $n\geq 2$, vector bundles over $\P^n_{\K}$ do not necessarily
split. Indeed, the tangent bundle is indecomposable.
In these cases, some sufficient conditions for vector bundles to
split have been established.
The following is one of such criterions, which we call
``Restriction criterion''.

\begin{theorem}[Horrocks]
\label{rest}
Let $\K$ be an algebraically closed field,
$n$ be an integer greater than or equal to 3, 
and let $E$ be a locally free sheaf on $\P^n_\K$ of $\rank\ r\ (\ge 1)$.
Then
$E$
splits
into a direct sum of line bundles if and only if
there exists a hyperplane $H \subset \P^n_\K$ such that
$E|_H$ splits into a direct sum of line bundles.
\end{theorem}
In other words,
the splitting of a vector bundle can be characterized
by using a hyperplane section.
However, vector bundles, or equivalently locally free sheaves, form a
small class among all coherent sheaves. There are some
important wider classes of coherent sheaves, e.g.,
reflexive sheaves or torsion free sheaves.
The purpose of this article is to generalize the 
``Restriction criterion'' to 
one for
reflexive sheaves, and we also show that it fails in the class of
torsion free sheaves.
Our main theorem is as follows.

 \begin{theorem}\label{main}
 Let $\K$ be an algebraically closed field,
 $n$ be an integer greater than or equal to 3, 
 and let $E$ be a reflexive sheaf on $\P^n_\K$ of $\rank\ r\ (\ge 1)$. Then
 $E$
 splits
 into a direct sum of line bundles if and only if
 there exists a hyperplane $H \subset \P^n_\K$ such that
  $E|_H$ splits into a direct sum of line bundles.
 \end{theorem}

 We give two proofs for Theorem \ref{main}. The first proof
 is basically parallel to that of Theorem \ref{rest}, in which 
 we also establish a general principle that the
 structure of a reflexive sheaf can be recovered from
 its hyperplane section (Theorem \ref{sp}).

 The second proof is based on a cohomological characterization
 for a coherent sheaf to be locally free. By using it,
 the proof is reduced to Theorem \ref{rest}.

 The organization of this paper is as follows.
 In \S\ref{pre}, we recall some basic results on reflexive sheaves
 from \cite{H2}.
 In \S\ref{pf}, we give the first proof of the main theorem.
 In \S\ref{app}, we give the second proof by using a cohomological
 characterization for a coherent sheaf to be locally free. 

 To each hyperplane arrangement in a vector space,
 we can associate a reflexive sheaf over the projective space.
 The splitting of this reflexive sheaf defines an important 
 class of arrangements, namely, free arrangements.
 As an application of our main theorem,
 we give a criterion for an arrangement
 to be free in \S\ref{arr}, which has been also obtained in
 \cite{Y}.

 \textbf{Acknowledgement.}
 The authors learned results of \S\ref{app} from Professor
 F.-O. Schreyer. They are grateful to him.
 The authors also thank to Takeshi Abe and Florin Ambro for many
 helpful comments and pointing out mistakes in our draft.
 The second author was supported by the JSPS Research Fellowship for
 Young Scientists.

 \section{Preliminaries}\label{pre}

 In this section, we fix the notation and prepare some results for the proof
 of Theorem
 {\rmfamily \ref{main}}. We use the terms ``vector bundle''
 and ``locally free sheaf''
 interchangeably. The term ``variety'' means a integral
 scheme of finite type over a field.
 Let $X$ be a smooth variety of dimension $n$ over a
 field $\K$,
 where $n \ge 1$ and $\K$ is an
 algebraically closed field.
 For a coherent sheaf $E$ on $X$ we denote by
 $\Sing (E)$ the non-free locus of $E$, i.e.,
 $\Sing(E):=\{ x \in X|
 E_x\ \mbox{is not a free}\
 \mathcal{O}_{x,X} \mbox{-module}\}$. The dual of a coherent sheaf
 $E$ (on $X$) is denoted by $E^*$.

 In this article, we employ homological algebra to investigate
 properties of
 a coherent sheaf on a smooth variety $X$.
 Let us
 review some definitions and results. For a coherent sheaf $E$ on
 $X$ over $\K$ and for a point $x \in X$
 (denoted by $\depth_{\mathcal{O}_X} (E_x))$ as the length of
 a maximal $E_x$-regular sequence in $\mathcal{M}_x$, where
 $\mathcal{M}_x$ is the unique maximal ideal of a local ring
 $\mathcal{O}_{x, X}$.
 Moreover, we define the projective dimension of an
 $\mathcal{O}_{x,X}$-module $E_x$ (denoted by
 $\pd_{\mathcal{O}_{x,X}} (E_x))$
 as the length of a minimal free resolution of $E_x$ as an
 $\mathcal{O}_{x,X}$-module.
 It is known that every module which is finitely generated over a regular
 local ring has finite projective
 dimension. These two quantities
 are related by the famous Auslander-Buchsbaum formula as follows.
 \[
 \depth_{\mathcal{O}_{x,X}}(E_x)+ \pd_{\mathcal{O}_{x,X}}
 (E_x)=\dim
 \mathcal{O}_{x,X}.
 \]
 Hence it follows easily that a coherent sheaf $E$ on $X$ is locally
 free if and only if
 $\depth_{\mathcal{O}_{x,X}} (E_x)= \dim \mathcal{O}_{x,X}$ for
 all
 $x \in X$. For details and proofs, see \cite{M}.
 The projective dimension can also be characterized as follows
 (for example, see \cite{OSS} Chapter II).

 \begin{lemma}
 \label{chara}
 Let $X$ be a smooth variety and
 $E$ be a coherent sheaf on $X$. Then $\pd_{\mathcal{O}_{x,X}}(E_x)\leq q$
 if
 and only if for
 all $i>q$ we have
 $$
 \mathcal{E}xt^i_{\mathcal{O}_X}(E, \mathcal{O}_X)_x=0.
 $$
 \end{lemma}
 In particular, $E$ is locally free if and only if
 $
 \mathcal{E}xt^i_{\mathcal{O}_X}(E, \mathcal{O}_X)=0
 $
 for all $i>0$.

 Next, let us review definitions and results on reflexive sheaves on
 $\P^n_\K$. Reflexive sheaves form a category between torsion free sheaves
 and
 vector bundles.
 \begin{define}
 We say a coherent sheaf $E$ on $\P^n_\K$ is reflexive if the canonical
 morphism
 $E \rightarrow E^{**}$ is an isomorphism.
 \end{define}
 In this article, we use the following results on reflexive sheaves. For the
 proofs and details, see \cite{H2}.
 \begin{prop}[\cite{H2}, Proposition 1.3]\label{depth}
 A coherent sheaf $E$ on $\P^n_\K$ is reflexive if and only if $E$ is
 torsion
 free and
 $\depth_{\mathcal{O}_{x,\P^n_\K}}(E_x) \ge 2$ for all points $x \in
 \P^n_\K$
 such that
 $\dim \mathcal{O}_{x,\P^n_\K} \ge 2$.
 \end{prop}

 \begin{cor}[\cite{H2}, Corollary 1.4]
 $\codim_{\P^n_\K} \Sing (E) \ge 3$ for a reflexive sheaf $E$ on $\P^n_\K$.
 \end{cor}

 \begin{prop}[\cite{H2}, Proposition 1.6]\label{123}
 %For a coherent sheaf $E$ on $\P^n_\K$, the followings are equivalent.
 %%%%  COMMENT  %%%%
 %%%% Oda-Tex claims ``the followings'' is a typical miss-used.
 %%%%  REPLY%%%%
 %%%% Thanks for your pointing out. I sometimes see this kind of miss-using.
 %Hence
 %%%% I also sometimes uses this.
 For a coherent sheaf $E$ on $\P^n_\K$, the following are equivalent.
 \begin{itemize}
 \item[1. ] $E$ is reflexive.
 \item[2. ] $E$ is torsion free and normal.
 \item[3. ] $E$ is torsion free and for each open set
 $U \subset \P^n_\K$ and each closed set $Z$ in $U$
 satistying $\codim_{U}(Z) \ge 2$,
 we have $E|_U \simeq j_* (E|_{U\setminus Z})$, where
 $j: U \setminus Z \rightarrow Z$ is an open immersion.
 \end{itemize}
 \end{prop}

 \section{The first proof of Theorem \ref{main}}\label{pf}
 Let us prove Theorem {\rmfamily \ref{main}}.
 It suffices to show the ``if'' part of the statement.
 %At first, let us assume that $\dim (\Sing(E)) \ge 1$. Then any hyperplane
 %%%%  COMMENT  %%%%
 %%%% ``At first'' is used to express ``initial state'' or
 %%%% ``initial condition'' like ``When t=0''.
 %%%% ``First'' means ``the first of some ordered process''.
 %%%%  REPLY%%%%
 %%%%  I understand. Thank you.
 First, let us assume that $\dim (\Sing(E)) \ge 1$. Then any hyperplane $H
 \subset \P^n_\K$ intersects 
  $\Sing(E)$. Take a point $x \in H \cap \Sing(E)
 \neq \emptyset$. Note that $\depth_{\mathcal{O}_{x, \P^n_\K}}
 (E_x) \le \dim \mathcal{O}_{x,\P^n_\K}-1$.
 Since the equation $h \in \mathcal{O}_{x,\P^n_\K}$
  which defines $H$ at $x$
 is a regular element for the reflexive $\mathcal{O}_{x,\P^n_\K}$-module
 $E_x$, it follows that $\depth_{\mathcal{O}_{x,H}} (E|_H)_x < \dim
 \mathcal{O}_{x,\P^n_\K}-1=
 \dim \mathcal{O}_{x,H}$. From Auslander-Buchsbaum formula,
 we conclude that $E|_H$ can not even be locally free.
 %nor
 %a direct sum of line bundles.
 Hence we may assume that $\dim (\Sing(E)) =0$.

 The next lemma is a generalization of Theorem 2.5 in
 \cite{H2}.

 \begin{lemma}\label{h1}
 Let $E$ be a reflexive sheaf on $\P^n_\K$ ($n\geq 3$)
 with $\dim (\Sing(E)) =0$.
 Suppose the restriction $E|_H$ to a hyperplane $H$ splits into a direct
 sum of
 line bundles. Then
 $$
 H^1(\P^n_\K,E(k))=0, \mbox{ for all } k\in\Z.
 $$
 \end{lemma}
 {\bf Proof of Lemma \ref{h1}}.
 %To show this lemma,
 We use the long exact sequence associated with the
 short exact sequence
 \[
 0 \rightarrow E(k-1) \rightarrow E(k) \rightarrow E(k)|_H \rightarrow 0.
 \]
 Because $E(k)|_H$ is a direct sum of line bundles, it follows that
 $H^1(H, E(k)|_H)=0$. So
 we have surjections
 \begin{equation}\label{surj}
 H^1(\P^n_\K, E(k-1)) \twoheadrightarrow H^1(\P^n_\K,E(k)),\ \forall k \in
 \Z.
 \end{equation}
 %Let us show that
 %$H^1(\P^n_\K,E(k))=0$
 %for all $k \in \Z$.
 %To see this, let us use the spectral sequence of
 %To prove the above lemma,
 To see that these cohomology groups are equal to zero, 
let us consider the spectral 
 sequence of
 local and global Ext functors:
 \[
 E_2^{p,q}=H^p(\P^n_\K,\mathcal{E}xt_{\P^n_\K}^q(E,\omega)) \Rightarrow
 E^{p+q}=\Ext_{\P^n_\K}^{p+q}(E,\omega)
 \]
 where $\omega$ is the dualizing sheaf of $\P^n_\K$. The assumption
 $\dim(\Sing (E))=0$ implies $\dim
 (\Supp(\mathcal{E}xt_{\P^n_\K}^q(E,\omega))) =0$ for all $q>0$. Thus it
 follows
 that $E_2^{p,q}=0$ unless
 $p =0$ or $q = 0$. Moreover,
 %since $E$ is reflexive and applying
 Proposition \ref{depth} implies
 $\depth_{\mathcal{O}_{x,\P^n_\K}}(E_x) \ge 2$. From
 Auslander-Buchsbaum formula, we have
 $\pd_{\mathcal{O}_{x,\P^n_\K}}E_x < n-1$ for all $x \in \P^n_\K$.
 It follows that
 $\mathcal{E}xt_{\P^n_\K}^q(E,\omega)=0$ for $\forall q\geq n-1$.
 Hence
 we have
 $E_2^{p,q}=0$ for $q \ge n-1$. Considering the convergence of this spectral
 sequence,
 we obtain the surjection
 \begin{equation}\label{surj2}
 H^{n-1}(\P^n_\K,\mathcal{H}om_{\P^n_\K}(E,\omega)) \simeq
 H^{n-1}(\P^n_\K,E^* \otimes \omega) \twoheadrightarrow
 \Ext^{n-1}_{\P^n_\K}(E,\omega).
 \end{equation}
 Since $\Ext^{n-1}_{\P^n_\K}(E(k),\omega)$ is the Serre dual to
 $H^1(\P^n_\K, E(k))$, they have the same dimension.
 From (\ref{surj2}),
 we have
 \begin{equation}\label{ineq}
 \dim H^1(\P^n_\K,E(k))\le
 \dim H^{n-1}(\P^n_\K, E^*(-k)\otimes \omega)
 \end{equation}
 for all $k \in \Z$. The right hand side of (\ref{ineq})
 vanishes for $k\ll 0$. Then together with the surjectivity (\ref{surj}),
 we conclude that
 $H^1(\P^n_\K,E(k))=0, \mbox{ for all } k\in\Z$.
 %This completes the proof of Lemma \ref{h1}.
 \owari

 %Combining this with the previous discussion, we can see that
 %$H^1(\P^n_\K,E(k))=0$ for all
 %$k \in \Z$.

 Now, let us put
 \[
 E|_H \simeq \ \oplus_{i=1}^r \mathcal{O}_H(a_i)
 \]
 %for integers $a_1 \ge a_2 \ge \cdots \ge a_r$
 %and for integers $a_1 \ge a_2 \ge \cdots \ge a_r$
 %%%%  COMMENT  %%%%
 %%%% Do we need the above assumption on degrees?
 %%%%  REPLY%%%%
 %%%% No, we do not. As you did, let us omit the assumption on degrees.
 and
 $F :=\oplus_{i=1}^r \mathcal{O}_{\P^n_\K}(a_i)$. Noting that
 $\Ext_{\P^n_\K}^1(F,E(-1)) \simeq H^1(\P^n_\K, E(-a_i-1))=0$,
  Theorem \ref{main} follows from the following theorem,
 which asserts that, roughly speaking, the structure of a reflexive sheaf
 can be recovered from its restriction to a hyperplane.

 \begin{theorem}\label{sp}
 Let $E$ and $F$ be
 reflexive sheaves on $\P^n_\K\ ( n \ge 2)$  and $H$ be a hyperplane in
 $\P^n_\K$. Suppose
 $E|_H \cong F|_H$ and $\Ext_{\P^n_\K}^1(F,E(-1))=0$. Then
 $E\cong F$.
 \end{theorem}
 {\bf Proof of Theorem \ref{sp}}.
 We want to extend the
 isomorphism
 $\varphi:F|_H \rightarrow E|_H$ to one over $\P^n_\K$. That is possible
 since there
 is an exact sequence
 \begin{eqnarray*}
 0 &\rightarrow& \Hom_{\P^n_\K}(F,E(-1)) \rightarrow \Hom_{\P^n_\K}(F,E)
 \rightarrow \Hom_{\P^n_\K}(F,E|_H)\\
 &\rightarrow& \Ext_{\P^n_\K}^1(F,E(-1)) =0,
 \end{eqnarray*}
 and every morphism $F|_H \rightarrow E|_H$ has a canonical extension to a
 morphism $F
 \rightarrow E|_H$. Let us fix 
 an extended morphism $f:F \rightarrow E$ which satisfies $f|_H =\varphi$.
 Now, let us
 consider the morphism $\det f : \det F \rightarrow \det E$. This
 is a monomorphism because $f$ is already a monomorphism.
 %Note that
 %the first Chern classes of $E$ and $F$ are the same.
 %This can be seen if we take the locally free resolution of $E$ and
%%% restrict
 %it
 %on any hyperplanes, then the restricted sequence is still exact
 %(Since $E$ is torsion free, there are no associated points. Hence
 %the exactness are kept when restricted according to Schlessinger's lemma).
 Since $E|_H \simeq F|_H$, ranks and first Chern classes of $E$ and $F$ are
 the same. Henceforth
  we can see that $\det f$ is a multiplication
 of some constant element in $\K$. Note that this constant is not zero. For
 $\det f$ is not zero on $H$. Thus at each point $x \in \P^n_\K \setminus
 (\Sing
 (E) \cup \Sing(F))$, the
 morphism $f_x$ is an isomorphism because at these points $f_x$
 are the endomorphism of a direct sum of local rings of the same rank.
 Since $\codim_{\P^n_\K}(\Sing(E) \cup \Sing(F)) > 2$ and both of
 $E$ and $F$ are reflexive, the third condition of Proposition \ref{123}
 implies
 that
 $f$ is also an isomorphism on $\P^n_\K$. \owari
 %Q.E.D.

 \begin{rem}
 In Theorem {\rmfamily \ref{main}}, we can not omit the assumption that
 $E$ is
 reflexive, i.e.,
 %e.g.,
 %the statement is not true for torsion free, but not reflexive sheaves.
 %%%%  ??
 ``Restriction criterion'' fails for torsion free sheaves.
 %To see this, we should
 For example, consider the ideal sheaf $I_p$ on $\P^3_{\K}$ which
 corresponds to a closed point $p \in \P^3_{\K}$. Note that
 $I_p$ is not reflexive. Indeed, let us put
 $U=\P^3_{\K} \setminus \{p\}$ and $j:U \rightarrow \P^3_{\K}$
 be an open immersion. It is easy to see that
 $I_p|_U \simeq \mathcal{O}_U$. If $I_p$ is reflexive, then according to
 Proposition \ref{123},
 $j_* (I_p|_U) \simeq I_p$ must hold. However, clearly this is not ture.
 Hence $I_p$ is not reflexive.
 Now,
 %the non-free locus of $I_p$ consists of only one
 %point $p$ (in particular,
 % $\dim (\Sing (E))=0)$.
 %If
 if we cut $I_p$ by a plane $H$ which does not contain $p$,
 then it is easily seen that
 $I_p|_H \simeq \mathcal{O}_H$. However, of course, $I_p$ is not a line
 bundle on $\P^3$.
 \end{rem}

 \section{The second proof}\label{app}
 Instead of Theorem {\rmfamily \ref{sp}}, we can use the following result,
 which is the generalization of the famous Horrocks' splitting criterion
 (For example,
 see \cite{OSS}). Combining this criterion with usual cohomological
 arguments and Lemma \ref{h1},
 we can give the second proof of Theorem {\rmfamily \ref{main}}. However, it
 seems
 that this theorem is not so familiar. Hence let us show the result with a
 complete proof. 

 \begin{theorem}\label{gspc}
 Let $\K$ be an algebraically closed field, $n$ be a integer greater 
than or equal to 2, 
  and let $E$ be a coherent sheaf on $\P^n_\K$. Then $E$ splits into a
 direct sum of line bundles if and only if $H^i(\P^n_\K,E(k))=0$ for all $k
 \in \Z,\ i=1,\cdots,n-1$ and $H^0(\P^n_{\K}, E(k))=0$ for all $k \ll 0$.
 \end{theorem}

 \begin{rem}
 Note that when $E$ is torsion free, then $H^0(\P^n_{\K}, E(k))=0$ for all
 $k \ll 0$. This follows from the fact that all torsion free sheaves 
 can be embedded into a direct sum of line bundles on $\P^n_{\K}$. So in the
 theorem, the condition $H^0(\P^n_{\K}, E(k))=0$
 is automatically satisfied for torsion free sheaves.
 \end{rem}
 When $E$ is a vector bundle, Theorem \rmfamily \ref{gspc} is just the
 splitting criterion of
 Horrocks. Thus for the proof of this theorem, it suffices to show the
 following lemma.

 \begin{lemma}\label{schreyer}
 Let $X$ be a nonsingular projective variety over an algebraically
 closed field $\K$ of dimension $n>1$, $L$ be an ample line bundle on $X$,
 and let $E$ be a
 coherent sheaf on $X$.Then
 $E$ is locally free if and only if $H^i(X, E(k))=0$ for all $k \ll 0$ and
 $i=0,1,\cdots,n-1$, where $E(k)=E\otimes L^k$. 
 \end{lemma}
 {\bf Proof of Lemma \ref{schreyer}}.
 From Serre duality, the ``only if'' part follows immediately. Let us show
 the
 ``if'' part of the statement.
 Recall that $E$ is locally free on $X$ if and only
 if
 $\mathcal{E}xt^i_X(E, \mathcal{O}_X)=0$
 for all $i >0$,
 %see Chapter II of [OSS].
 see \S\ref{pre}.
 %(For example, see \cite{OSS}).
 Consider the spectral sequence
 \[
 E_2^{p,q}(k)=H^p(X,\mathcal{E}xt^q_X(E(k),\omega)) \Rightarrow
 E^{p+q}(k)=\Ext^{p+q}_X
 (E(k),\omega),
 \]
 where $k \in \Z$ and $\omega$ is the dualizing sheaf on $X$.
 % and $\omega \simeq \mathcal{O}_{\P^n_\K}(-n-1)$ is a
 %dualizing sheaf on $\P^n_\K$.
 %%%%% We once defined $\omega$.
 By Serre
 duality, $H^i(X,E(k))^* \simeq \Ext^{n-i}_X(E(k), \omega)$ for
 $i=0,1,\cdots, n$.
  So for each $i > 0,\ E^i(k)= \Ext^i_X(E(k),\omega) = 0$ for sufficiently
 small
  $k \in \Z$. Now let us
 assume that there exists an integer $i>0$ such that
 $\mathcal{E}xt^i_X(E,\mathcal{O}_X) \neq 0$,
  and we show that this leads to a contradiction.
 It is easy to see that
 %there exists an integer $k_0$ such that
 \[
 E_2^{0,i}(k)=H^0(X, \mathcal{E}xt_X^i(E, \omega) \otimes
 \mathcal{O}_X(-k)) \neq 0,\ \mbox{for }\forall k\ll 0.
 \]
 %for all $k<k_0$.
 %Since $\mathcal{E}xt^i_{\P^n_\K} (E,
 %\mathcal{O}_{\P^n_\K}) \neq 0$ and
 %$E_2^{0,i}(k)$ is a subspace of $E_2^{0,i}(k-1)$, there exists an integer
 %$k_0$ such that for all $k < k_0,\ E_2^{0,i}(k) \neq 0$.
 %Next, let us
 %consider the differential map
 %$d_2^{0,i}(k): E_2^{0,i}(k) \rightarrow E_2^{2,i-1}(k)$. Since
 %$E_2^{2,i-1}(k)=
 %H^2(\P^n_\K, \mathcal{E}xt^{i-1}_{\P^n_\K}(E(k), \omega))$, there exists
% an
 %integer $k_2$ such that for all $k < k_2,
 %\ E_2^{2,i-1}(k) =0$. So for sufficiently small $k$, we can see that
 %$E_2^{0,i}(k)=E_3^{0,i}(k)$.
 %Continuing this process, we can find an integer $k'$ such that for all
 %$k < k',\ E_2^{0,i}(k) = E_{\infty}^{0,i}(k) \neq 0$.
 On the other hand, for
 $p>0$,
 $$
 E_2^{p,q}(k)=H^p(X,\mathcal{E}xt^q_X(E,\omega)\otimes
 \mathcal{O}_X(-k))
 =0,\ \mbox{for }\forall k\ll 0.
 $$
 From the definition of spectral sequence,
 $$
 \Ext ^i_X(E(k), \omega)=E_2^{0,i}(k)\neq 0,
 $$
 for $\forall k\ll 0$.
 This contradicts the
 assumption that
 for each $i>0,\ E^i(k) = 0$ for sufficiently small $k \in \Z$. Hence we
 can see that
 $\mathcal{E}xt^i_X(E, \mathcal{O}_X)=0$ for all
 $i>0$, so $E$ is a locally free sheaf. \owari
 %%Q.E.D.

 \section{Application to hyperplane arrangements}\label{arr}

 In this section, we describe an application of our main theorem
 to the theory of hyperplane arrangements. As mentioned in
 \S \ref{intro}, each hyperplane arrangement determines a
 reflexive sheaf. We start with a 
 more general setting. To every divisor $D$ in a complex manifold $M$
 we can associate a reflexive sheaf as follows.

 \begin{define}
 \label{log}
 A vector field $\delta$ on an open set $U\subset M$ is
 said to be logarithmic tangent to $D$ if for a local defining
 equation $h$ of $D\cap U$ on $U$, $\delta h\in (h)$.
 The sheaf associated with logarithmic vector fields
 is denoted by $\Der_M(-\log D)$.
 \end{define}
 In the definition above, a vector field $\delta$ is identified with
 a derivation $\delta:\mathcal{O}_M\longrightarrow \mathcal{O}_M$, and
 $\Der_M(-\log D)$ can be considered as a subsheaf of the tangent sheaf.
 The sheaf of logarithmic vector fields $\Der_M(-\log D)$ is
 not necessarily locally free, but in \cite{Slog}, K. Saito proved
 the following.

 \begin{theorem}[\cite{Slog}]
 $\Der_M(-\log D)$ is a reflexive sheaf.
 \end{theorem}
 From now on,
 %%%% A2 The change here is just only my favor.
  we restrict ourselves to the case where $D$
 is a hyperplane arrangement.

 Let $V$ be an $\ell$-dimensional linear space over
 $\K$ and
 $S:=\K[V^*]$ be the algebra of polynomial functions on $V$ that is
 naturally isomorphic to $\K[z_1, z_2, \cdots, z_\ell]$ for
 any choice of basis $(z_1, \cdots, z_\ell)$ of $V^*$.

 A (central) hyperplane arrangement $\mathcal{A}$ is a finite
 collection of codimension one
 linear subspaces in $V$.
 For each hyperplane $H$ of $\mathcal{A}$, fix
 a nonzero linear form $\alpha_H\in V^*$ vanishing on $H$ and
 put $Q:=\prod_{H\in\mathcal{A}}\alpha_H$.

 The characteristic polynomial of $\mathcal{A}$ is defined as
 $$
 \chi(\mathcal{A}, t)=\sum_{X\in L_\mathcal{A}}\mu(X)t^{\dim X},
 $$
 where $L_\mathcal{A}$ is a lattice which consists of the intersections of
 elements of $\mathcal{A}$, ordered by reverse inclusion,
 $\hat{0}:=V$ is the unique minimal element of $L_\mathcal{A}$ and
 $\mu:L_\mathcal{A}\longrightarrow\Z$ is the M\"obius function
 defined as follows:
 \begin{eqnarray*}
 \mu(\hat{0})&=&1,\\
 \mu(X)&=&-\sum_{Y<X}\mu(Y),\ \mbox{if}\ \hat{0}<X.
 \end{eqnarray*}
 The characteristic polynomial is one of the most important
 concepts in the theory of hyperplane arrangements. Actually
 there are a lot of combinatorial or geometric interpretations
 of characteristic polynimial.
 For details, see \cite{OT}.

 Denote by $\Der_V:=\K[V^*]\otimes V$ the $S$-module of
 all polynomial vector fields on $V$. The following definition
 was given by G. Ziegler.

 \begin{define}[\cite{Z}]
 For a given arrangement $\mathcal{A}$ and a map
 $m:\mathcal{A}\longrightarrow\Z_{\geq 0}$, we define modules
 of logarithmic vector fields with multiplicity $m$ by
 $$
 D(\mathcal{A}, m)=\{\delta\in\Der_V\ |\ \delta \alpha_H \in S
 \alpha^{m(H)},\
 \forall H\in\mathcal{A}\}
 $$
 When the multiplicity $m$ is the constant map $\underbar{{\rm 1}}(H)\equiv
 1\
 (\forall H\in\mathcal{A})$, $D(\mathcal{A}, \underbar{{\rm 1}})$ is simply
 denoted by $D(\mathcal{A})$.
 \end{define}
 It is known that the graded $S$-module $D(\mathcal{A}, m)$ is a reflexive
 module of rank $l= \dim V$.
 %Next, let us introduce the concept of free arranegments, which plays an
 %important role
 %in the study of arrangements.

 \begin{define}
 \begin{itemize}
 \item[(1)]
 An arrangement with a multiplicity $(\mathcal{A}, m)$ is called
 free with exponents $(e_1, \cdots, e_\ell)$ if $D(\mathcal{A}, m)$ is
 a free $S$-module, with a homogeneous basis
 $\delta_1, \cdots, \delta_\ell$ such that
 $$
 \deg \delta_i=e_i.
 $$
 Note that a vector field
 $$
 \delta=\sum_i f_i\frac{\partial}{\partial x_i}
 $$
 is said to be homogeneous if coefficients $f_1, \cdots, f_\ell$ are
 all homogeneous with the same degree
 and put $\deg \delta :=\deg f_i$.
 \item[(2)]
 An arrangement $\mathcal{A}$ is called free if
 $(\mathcal{A}, \underbar{{\rm 1}})$ is free, i.e.,
 $D(\mathcal{A})$ is a free $S$-module.
 \end{itemize}
 \end{define}
 Since $D(\mathcal{A})$ contains the Euler vector field
 $\theta_E:=\sum_{i=1}^\ell x_i\frac{\partial}{\partial x_i}$,
 the exponents $(e_1, \cdots, e_\ell)$ of a free arrangement
 $\mathcal{A}$ contains $1$.
 H. Terao proved that the freeness of $\mathcal{A}$
 implies a remarkable behavior of the characteristic polynomial.

 \begin{theorem}[\cite{Tfact}]
 \label{factor}
 Suppose $\mathcal{A}$ is a free arrangement with the exponents
 $(e_1, \cdots, e_\ell)$, then
 $$
 \chi(\mathcal{A}, t)=\prod_{i=1}^\ell(t-e_i).
 $$
 \end{theorem}

 As we will see later, in Corollary \ref{free},
 the freeness is equivalent to the splitting of a reflexive sheaf, and
 exponents are corresponding to the splitting type. On the other hand,
 the left hand side of the Theorem \ref{factor}
 is obtained from the intersection poset, thus 
 determined by the combinatorial structure.
 This theorem connects two regions in mathematics: 
 combinatorics of arrangements and geometry of reflexive sheaves.
 It enables us to study combinatorics of arrangements via a 
geometric method. 
 For example, in \cite{Y} characteristic polynomials for some
 arrangements are computed by using this interpretation.

 In \cite{Z}, Ziegler studied the relation between
 the freeness and the freeness with a multiplicity.
 Fixing a hyperplane $H_0\in\mathcal{A}$, let us define
 an arrangement
 $$
 \mathcal{A}^{H_0}:=\{ H_0\cap K\ |\ K\in\mathcal{A},\ K\neq H_0\},
 $$
 over $H$ and the natural multiplicity
 $$
 \underline{m}(X):=\sharp\{ K\in\mathcal{A}\ |\ K\cap H_0=X\}
 $$
 for $X\in\mathcal{A}^{H_0}$.

 \begin{theorem}[\cite{Z}]
 \label{thm:zie}
 If $\mathcal{A}$ is a free arrangement with
 exponents $(1, e_2, \cdots, e_\ell)$, then the restricted
 arrangement with natural multiplicity
 $(\mathcal{A}^{H_0}, \underbar{m})$ is also free
 with exponents $(e_2, \cdots, e_\ell)$.
 \end{theorem}

 More precisely,
 let $\alpha=\alpha_{H_0}$ be a defining equation of $H_0$ and define
 $$
 D_0(\mathcal{A}):=\{ \delta\in D(\mathcal{A})\ |\
 \delta\alpha=0\}.
 $$
 It is easily seen that $D(\mathcal{A})$ has a direct sum
 decomposition into graded $S$-modules
 $$
 D(\mathcal{A})=S\cdot \theta_E\oplus D_0(\mathcal{A}).
 $$
 Ziegler proved that if $\delta_1=\theta_E, \delta_2, \cdots, \delta_\ell$
 is a basis of $D(\mathcal{A})$ with
 $\delta_2, \cdots, \delta_\ell \in D_0(\mathcal{A})$, then
 $\delta_2|_{H_0}, \cdots, \delta_\ell|_{H_0}$ form a basis
 of $D(\mathcal{A}^{H_0}, \underbar{m})$.

 Recall that a graded $S$-module $M=\oplus_{k\in\Z}M_k$ determines
 a coherent sheaf $\tilde{M}$ over $\P^{\ell-1}=\Proj S$. Conversely
 for any coherent sheaf $\mathcal{F}$ over $\P^{\ell-1}$,
 $\Gamma_*(\mathcal{F}):=\bigoplus_{k\in\Z}\Gamma(\P^{\ell-1},
 \mathcal{F}(k))$
 defines the graded $S$-module associated with $\mathcal{F}$.
 We have the natural $S$-homomorphism
 $\alpha:M\rightarrow \Gamma_*(\tilde{M})$, which is neither
 injective nor surjective in general. In the case of $M=D(\mathcal{A})$,
 however,
 we have the following lemma.

 \begin{lemma}
 \label{graded}
 $
 \alpha: D(\mathcal{A})\stackrel{\cong}{\longrightarrow}
 \Gamma_*\left(\P^{\ell-1}, \widetilde{D(\mathcal{A})}\right)$
 is isomorphic.
 \end{lemma}
 {\bf Proof of Lemma \ref{graded}}. We prove the surjectivity.
 Since $\bigcup_{i=1}^\ell D(z_i)=\P^{\ell-1}$, any element
 in $\Gamma(\P^{\ell-1}, \widetilde{D(\mathcal{A})}(k))$ can be
 expressed as
 $$
 \delta=
 \frac{\delta_1}{z_1^{d_1}}=
 \frac{\delta_2}{z_2^{d_2}}=\cdots =
 \frac{\delta_\ell}{z_\ell^{d_\ell}},
 $$
 where $\delta_i\in D(\mathcal{A})_{d_i+k}$. From the facts that
 $\delta_i$ is an element of a $S$-free module $\Der_V$ and
 $S$ is UFD,
 it is easily seen that $\delta$ is also a polynomial
 vector field, so contained in $\Der_V$.
 Let $\alpha_H$ be a defining linear form
 of $H\in\mathcal{A}$, and we may choose $i$ such that
 $\alpha_H$ and $z_i$
 are linearly independent. Then the right hand side of
 $$
 z_i^{d_i}\cdot \delta\alpha_H=\delta_i\alpha_H
 $$
 is divisible by $\alpha_H$, so is the left.
 Hence $\delta\alpha_H$ is also divisible by
 $\alpha_H$, and we can conclude that $\delta\in D(\mathcal{A})$.
 \owari

 The above lemma enable us to connect freeness and splitting.

 \begin{cor}
 \label{free}
 $\mathcal{A}$ is free with exponents $(e_1, \cdots, e_\ell)$ if
 and only if
 $$
 \widetilde{D(\mathcal{A})}=
 \mathcal{O}_{\P^{\ell-1}}(-e_1)\oplus\cdots\oplus\mathcal{O}_{\P^{\ell-1}}(-e_\ell)
 $$
 \end{cor}

 Now, the following theorem, which has been proved and
 played an important role in the proof of Edelman and Reiner conjecture in
 \cite{Y},
 is naturally proved from Theorem \ref{main}.

 \begin{theorem}[\cite{Y}]
 \label{thm:y}
 $\mathcal{A}$ is free if and only if
 there exists a hyperplane $H_0\in\mathcal{A}$ such that
 \begin{itemize}
 \item[(a)] $(\mathcal{A}^{H_0}, \underbar{m})$ is free, and
 \item[(b)] $\mathcal{A}_x:=\{H\in\mathcal{A}\ |\ H\ni x\}$ is
 free for all $x\in H_0\setminus \{0\}$.
 \end{itemize}
 \end{theorem}
 {\bf Proof of Theorem \ref{thm:y}}.
 Let us denote by $\P(V)$ the
 projective space of one-dimensional subspaces in a vector space $V$.
 Recall that $D_0(\mathcal{A})$ is a graded reflexive $S$-module.
 So it determines a reflexive sheaf $\widetilde{D_0(\mathcal{A})}$ over
 $\P(V)$. As is mentioned in \cite{MS}, the local structure of
 $\widetilde{D_0(\mathcal{A})}$ is determined by the local structure
 of $\mathcal{A}$, i.e.,
 $$
 \widetilde{D_0(\mathcal{A})}_{\bar{x}}=
 \widetilde{D_0(\mathcal{A}_x)}_{\bar{x}},
 $$
 for $\bar{x}\in\P(V)$. Using Theorem \ref{thm:zie} locally,
 condition (b) in Theorem \ref{thm:y} implies that
 $$
 \widetilde{D_0(\mathcal{A})}_{\bar{x}}|_{\P(H_0)}
 =
 \widetilde{D(\mathcal{A}^{H_0}, \underbar{m})}_{\bar{x}}.
 $$
 Now condition (a) in Theorem \ref{thm:zie} means that
 $\widetilde{D_0(\mathcal{A})}|_{\P(H_0)}$ splits into
 a direct sum of line bundles.
 From Theorem \ref{main}, we may conclude that
 $\widetilde{D_0(\mathcal{A})}$ is also splitting.
 Hence
 $$
 \bigoplus_{k\in\Z}\Gamma\left(\P(V), \widetilde{D_0(\mathcal{A})}(k)\right)
 =D_0(\mathcal{A})
 $$
 is a free module over $S$.
 Thus $\mathcal{A}$ is a free arrangement.

 \owari

 %owari

 \noindent
 Takuro Abe\\
 Department of Mathematics, \\
 Kyoto University, \\
 Kyoto 606-8502, \\
 Japan, \\
 abetaku@math.kyoto-u.ac.jp

 \vspace{3mm}

 \noindent
 Masahiko Yoshinaga\\
 Research Institute for Mathematical Sciences, \\
 Kyoto University, \\
 Kyoto 606-8502, \\
 Japan\\
 yosinaga@kurims.kyoto-u.ac.jp
\end{document}